\documentclass[11pt]{article}
\usepackage[utf8]{inputenc}
\usepackage{amssymb,amsmath}
\usepackage{graphicx,units,mathrsfs,bm}
\usepackage{subfigure}
\usepackage[sort&compress,numbers]{natbib} 
\usepackage{hyperref,xcolor}
\usepackage[paperwidth=8.5in, paperheight=11in, portrait, top=1in, bottom=1.in, left=1.15in, right=1.15in]{geometry}
\usepackage{tikz}
\usetikzlibrary{calc,decorations.markings,positioning}

\newcommand\btd{\raise 2pt \hbox{$\hat\bigtriangledown$}\hskip 1.5pt}
\newcommand\bt{\raise 2pt \hbox{$\bigtriangledown$}\hskip 1.5pt}

\def\osp{{\mathfrak{osp}}}

\def\os{{\mathfrak{o}}}
\def\su{{\mathfrak{su}}}

\def\mD{{\mathcal{D}}}
\def\mA{{\mathcal{A}}}
\def\mP{{\mathcal{P}}}
\def\mQ{{\mathcal{Q}}}
\def\mK{{\mathcal{K}}}
\def\mR{{\mathcal{R}}}

\def\cliff{{C\!\!\ell}}

\newcommand{\ket}[1]{|#1\rangle}

\numberwithin{equation}{section}

\begin{document}
\title{
\bf A Howe correspondence for the algebra of the $\mathfrak{osp}(1|2)$ Clebsch-Gordan coefficients}
\author{
Julien Gaboriaud\footnote{E-mail: julien.gaboriaud@umontreal.ca}~, Luc Vinet\footnote{E-mail: vinet@CRM.UMontreal.CA} \\
\small~Centre de Recherches Math\'ematiques, Universit\'e de Montr\'eal, \\
\small~P.O. Box 6128, Centre-ville Station, Montr\'eal (Qu\'ebec), H3C 3J7, Canada.\\[.9em]
}
\date{\today}
\maketitle

\hrule
\begin{abstract}\noindent
Two descriptions of the dual $-1$ Hahn algebra are presented and shown to be related under Howe duality.
The dual pair involved is formed by the Lie algebra $\mathfrak{o}(4)$ and the Lie superalgebra $\mathfrak{osp}(1|2)$.\\[1em]
\end{abstract} 
\hrule


\section{Introduction}
The dual $-1$ Hahn algebra \cite{Genest2013} captures the bispectrality properties of the orthogonal polynomials bearing the same name. These polynonomials were first obtained \cite{Tsujimoto2013} as a $q=-1$ limit of the dual $q$-Hahn polynomials and were shown to essentially define the Clebsch-Gordan coefficients of the Lie superalgebra $\mathfrak{osp}(1|2)$. This paper provides two related pictures of the (centrally extended) dual $-1$ Hahn algebra and explains their connection on the basis of Howe duality.

The dual $-1$ Hahn algebra is an example of the quadratic algebras of Askey-Wilson type that are realized by the recurrence and the difference/differential equation operators of the hypergeometric orthogonal polynomials of the Askey tableaux that correspond to $q=1$, $q$ generic \cite{Koekoek2010}  but also to $q=-1$ in the case of the so-called Bannai-Ito scheme. Sitting at the top of each of these are respectively, the Racah, the Askey-Wilson and the Bannai-Ito \cite{Tsujimoto2012} polynomials. It is known that these orthogonal polynomials are, in the order in which they are listed, the Racah coefficients of the Lie algebra $\mathfrak{su}(2)$ or $\mathfrak{su}(1,1)$, of the quantum algebra $U_q(\mathfrak{sl}(2))$ and of the superalgebra $\mathfrak{osp}(1|2)$. This hints to the proven fact that the algebras associated to each of these families of polynomials and called by their names are realized as centralizers of the diagonal action of either $\mathfrak{su}(1,1)$, $U_q(\mathfrak{sl}(2))$ or $\mathfrak{osp}(1|2)$ on their three-fold product, an observation which is paving the way to higher rank extensions \cite{DeBie2016a,DeBie2017a,Post2017,DeBie2018}.

This feature explains in part why these algebras of Askey-Wilson type have become preeminent in a number of areas in mathematics and physics such as representation theory  
\cite{Crampe2019,Crampe2019d,Crampe2020}, combinatorics \cite{Terwilliger2014}, knot theory \cite{Bullock1999} and integrable models \cite{Genest2014,DeBie2015,Baseilhac2005}. There is hence much interest in deepening their understanding. In this respect, the Racah, Bannai-Ito and Askey-Wilson algebras have been given complementary descriptions \cite{Gaboriaud2018,Gaboriaud2018a,Gaboriaud2018b,Frappat2019} as commutants of maximal Abelian subalgebras in the universal algebra $U(\mathfrak{o}(n))$ of the orthogonal algebra $\mathfrak{o}(n)$ and its non standard $q$-deformation denoted $\mathfrak{o}_q(n)$ or $U_q'(\mathfrak{o}(n))$. Furthermore, it has been observed that these alternative presentations are in Howe duality with the centralizer ones. For a recent review see \cite{Gaboriaud2019}.

The dual Hahn polynomials, a limit of the Racah polynomials, enter in the Clebsch-Gordan coefficients of $\mathfrak{su}(1,1)$. Their algebra which is isomorphic to the Higgs algebra \cite{Higgs1979} arises in many contexts (see references in \cite{Frappat2019a}) and in particular as a truncation of the reflection algebra \cite{Crampe2019}. This dual Hahn algebra also admits two presentations which are related through Howe duality based on the pair $(\mathfrak{o}(4), \mathfrak{su}(1,1))$ \cite{Frappat2019a}: in the first one, it is the commutant of $\mathfrak{o}(2) \oplus \mathfrak{o}(2)$ in the oscillator representation of $U(\mathfrak{u}(4))$ and in the second it is embedded in $U(\mathfrak{su}(1,1))\otimes U(\mathfrak{su}(1,1))$. A  $q$-deformation of this analysis was performed \cite{Frappat2019} to define a $q$-analog of the Higgs algebra and provide a two-fold description of the dual $q$-Hahn algebra with Howe duality resting in this case on the pair $(\mathfrak{o}_q(4), U_q(\mathfrak{su}(1,1))$. Hence at the second level of the (discrete) Askey tableaux, this leaves the case $q = -1$, that is the algebra of the dual $-1$ Hahn polynomials, as the only one for which a description in the framework of Howe duality has not been given. The purpose of the present paper is to fill this gap.

The dual $-1$ Hahn algebra has been characterized in \cite{Genest2013}; it has been shown to arise in the studies \cite{Genest2012,Genest2014a} of the Dunkl oscillator in the plane and quite recently as the symmetry algebra of a superintegrable two-dimensional singular oscillator with internal degrees of freedom \cite{Bernard2020}. The dual $-1$ Hahn polynomials form the Clebsch-Gordan coefficients of $\mathfrak{osp}(1|2)$. This superalgebra will be one element of the dual pair at play, the other will be $\mathfrak{o}(4)$. It will be seen that the dual $-1$ Hahn algebra is the commutant of $\mathfrak{o}(2) \oplus \mathfrak{o}(2)$ in a spinorial representation of $\mathfrak{o}(4)$ given in terms of Bosonic and Fermionic oscillators. This picture will be shown to be in a duality relation with the embedding of the dual $-1$ Hahn algebra in $U(\mathfrak{osp}(1|2)\otimes U(\mathfrak{osp}(1|2)$ given in \cite{Genest2013}.

The paper will unfold as follows. Section \ref{sec:dualscwdunkl} will recall the definition of the dual $-1$ Hahn algebra and introduce as well the Schwinger-Dunkl algebra $\mathfrak{sd}(2)$ which was identified \cite{Genest2012,Genest2014a} as the symmetry algebra of the Dunkl oscillator in two dimensions. The nomenclature comes from the fact that $\mathfrak{sd}(2)$ is obtained when the raising operators in the Schwinger construction of $\mathfrak{su}(2)$ are replaced by creation and annihilation operators involving Dunkl operators instead of ordinary derivatives. Section \ref{sec:dualcommutant} will confirm that the dual $-1$ Hahn algebra can be realized as a commutant in the fashion described above. The commutant will initially be identified as $\mathfrak{sd}(2)$ and this is why the relation with the dual $-1$ Hahn algebra will have been established in the preceding section. Section \ref{sec:clebschgordan} will recall how the dual $-1$ Hahn algebra is embedded in  $U(\mathfrak{osp}(1|2) \otimes U(\mathfrak{osp}(1|2)$ in light of the fact that this dual $-1$ Hahn algebra characterizes the Clebsch-Gordan coefficients of $\mathfrak{osp}(1|2)$. That Howe duality connects the commutant and the embedding presentations of the dual $-1$ Hahn algebra will be the subject of Section \ref{sec:howeduality} and concluding remarks will form Section \ref{sec:conclusion}.

\section{The dual $-1$ Hahn algebra and the Schwinger-Dunkl algebra}\label{sec:dualscwdunkl}
We first introduce the two algebras of interest and show how they are closely related to each other.

\subsection{The dual $-1$ Hahn algebra}
The dual $-1$ Hahn algebra is defined \cite{Genest2013} by the generators $\mathbf{P}$, $K_1$, $K_2$, $K_3$ and the relations
\begin{align}\label{eq:dualm1Hahn}
\begin{aligned}
 &[K_1,K_2]=K_3,\qquad\qquad [K_1,K_3]=K_2+\nu \mathbf{P}+\tfrac{1}{2},\\
 &[K_2,K_3]=4K_1(1+\nu \mathbf{P})-2\nu K_3\mathbf{P}+\sigma \mathbf{P}+\rho,\\
 &[K_1,\mathbf{P}]=0,\qquad\qquad \{K_2,\mathbf{P}\}=-\mathbf{P}-2\nu,\qquad\qquad \{K_3,\mathbf{P}\}=0,
\end{aligned}
\end{align}
where $\nu$, $\sigma$, $\rho$ are structure constants. By promoting the structure constants to central elements, one obtains what will be referred to as the centrally extended dual $-1$ Hahn algebra. Note that $\rho$ can be reabsorbed in the generator $K_1$.

The dual $-1$ Hahn algebra captures the bispectral properties of the polynomials with the same name \cite{Tsujimoto2013} through a realization where $K_1$ is the $5$-diagonal difference of which they are eigenfunctions, $K_2$ is associated to the recurrence relation and $P$ is the parity or reflection operator.
%

\subsection{The Schwinger-Dunkl algebra}
The Schwinger-Dunkl algebra is the symmetry algebra of a two-dimensional isotropic Dunkl oscillator in the plane. This system is described by the Hamiltonian
\begin{align}
 H_{12}=-\tfrac{1}{2}\left[(\mD_{x_1}^{\mu_1})^{2}+(\mD_{x_2}^{\mu_2})^{2}\right]+\tfrac{1}{2}\left[{x_1}^{2}+{x_2}^{2}\right]=H_{1}+H_{2},\qquad \mD_{x_i}^{\mu_{i}}=\partial_{x_{i}}+\frac{\mu_{i}}{x_{i}}(I-R_{i}),
\end{align}
where $I$ is the identity operator and the $R_{i}$, $i\in\{1,2\}$ are the reflection operators
\begin{align}
 R_{1}f(x_1,x_2)=f(-x_1,x_2),\qquad R_{2}f(x_1,x_2)=f(x_1,-x_2).
\end{align}
The symmetry algebra of this system is obtained through the Schwinger construction. Form
\begin{align}
 \mathbf{a}_{i}^{\dagger}=\frac{x_i-\mD_{x_i}^{\mu_{i}}}{\sqrt{2}},\qquad \mathbf{a}_{i}^{\phantom{\dagger}}=\frac{x_i+\mD_{x_i}^{\mu_{i}}}{\sqrt{2}},
\end{align}
the parabosonic creation and annihilation operators, whose commutation relations are
\begin{align}
 [\mathbf{a}^{\phantom{\dagger}}_{i},\mathbf{a}_{j}^{\dagger}]=(I+2\mu_{i}R_{i})\delta_{ij}.
\end{align}
Then the three quantities
\begin{align}
 J_1=\tfrac{1}{2}\left(\mathbf{a}_1^{\dagger}\mathbf{a}_2^{\phantom{\dagger}}+\mathbf{a}_1^{\phantom{\dagger}}\mathbf{a}_{2}^{\dagger}\right),\qquad J_2=\tfrac{1}{2i}\left(\mathbf{a}_1^{\dagger}\mathbf{a}_2^{\phantom{\dagger}}-\mathbf{a}_1^{\phantom{\dagger}}\mathbf{a}_{2}^{\dagger}\right),\qquad J_3=\tfrac{1}{2}\left(H_1-H_2\right),
\end{align}
are symmetries of the Hamiltonian $H_{12}$, along with $R_1$, $R_2$. These elements $J_1$, $J_2$, $J_3$, $R_1$, $R_2$ obey the following commutation relations, which we will refer to as the relations of the Schwinger-Dunkl algebra $\mathfrak{sd}(2)$:
\begin{align}
\begin{aligned}{}
 &[J_2,J_3]=iJ_1,\qquad\qquad [J_3,J_1]=iJ_2,\\
 &[J_1,J_2]=i\left(J_3(1+\mu_1R_1+\mu_2R_2)-\tfrac{1}{2}H_{12}(\mu_1R_1-\mu_2R_2)\right),\\
 &\{J_1,R_\alpha\}=0,\qquad\qquad \{J_2,R_\alpha\}=0,\qquad\qquad [J_3,R_\alpha]=0,\qquad\qquad \alpha=1,2.
\end{aligned}
\end{align}

\subsection{Connection between the two algebras}
Starting from the generators $\mathbf{P}$, $K_1$, $K_2$, $K_3$, write
\begin{align}
 j_1=\tfrac{i}{2}K_3,\qquad j_2=-\tfrac{1}{2}(K_2+\nu \mathbf{P}+\tfrac{1}{2}),\qquad j_3=-K_1 -\frac{\rho}{4},
\end{align}
the dual $-1$ Hahn algebra relations now take the form
\begin{align}
\begin{aligned}{}
 &[j_2,j_3]=ij_1,\qquad\qquad [j_3,j_1]=ij_2,\\
 &[j_1,j_2]=i\left(j_3(1+\nu \mathbf{P})+\tfrac{1}{4}(\nu\rho-\sigma) \mathbf{P}\right),\\
 &\{j_1,\mathbf{P}\}=0,\qquad\qquad \{j_2,\mathbf{P}\}=0,\qquad\qquad [j_3,\mathbf{P}]=0.
\end{aligned}
\end{align}
One then sees that the dual $-1$ Hahn algebra is indeed similar to the $\mathfrak{sd}(2)$ algebra. The difference is that $\mathfrak{sd}(2)$ has $2$ reflection-type operators, $R_1$ and $R_2$, whilst the dual $-1$ Hahn algebra only has a single one, $\mathbf{P}$. It will thus prove more useful to work with $R_1$ and $R_{12}=R_1R_2$ as the latter commutes with everything and can be viewed as a central element in $\mathfrak{sd}(2)$.

Then, the two algebras obey the same relations upon identifying
\begin{align}
\begin{aligned}{}
 \mathbf{P}&=R_1,\\
 \nu&=\mu_1+\mu_2 R_{12},\\
 \rho&=2H_{12},\\
 \sigma&=2\mu_1\,\rho,
\end{aligned}
\end{align}
where we recall that $H_{12}$ and $R_{12}$ are both central elements.

Hence, the Schwinger-Dunkl algebra $\mathfrak{sd}(2)$ is essentially the centrally extended dual $-1$ Hahn algebra. In the remainder of the paper, we will encounter instances of algebras presented in the form of this $\mathfrak{sd}(2)$ algebra.

\section{The dual $-1$ Hahn algebra as a commutant}\label{sec:dualcommutant}
In this Section we will obtain the dual $-1$ Hahn algebra as the commutant of $\os(2)\oplus\os(2)$ in a spinorial realization.

\subsection{The model}
Consider the Hamiltonian
\begin{align}
 H=\frac{1}{2}\sum_{i=1}^{4}\{a_i^{\dagger},a_i^{\phantom{\dagger}}\},
\end{align}
built from the standard Bosonic raising $a_i^{\dagger}$ and lowering $a_i$ operators obeying
\begin{align}
\begin{aligned}{}
 [a_i^{\phantom{\dagger}},a_j^{\dagger}]&=\delta_{ij},
\end{aligned}\qquad
\begin{aligned}{}
 [a_i,a_j]&=[a_i^{\dagger},a_j^{\dagger}]=0.
\end{aligned}
\end{align}
We also introduce Fermionic raising $b_i^{\dagger}$ and lowering $b_i$ operators obeying
\begin{align}
\begin{aligned}{}
 \{b_i^{\phantom{\dagger}},b_j^{\dagger}\}&=\delta_{ij},
\end{aligned}\qquad
\begin{aligned}{}
 \{b_i,b_j\}&=\{b_i^{\dagger},b_j^{\dagger}\}=0.
\end{aligned}
\end{align}
Above and below $i, j=1, 2, 3, 4$. The Bosonic and Fermionic operators mutually commute with each other. One sees that the combinations 
\begin{align}
 \gamma_i=b_i^{\phantom{\dagger}}+b_i^{\dagger}
\end{align}
obey
\begin{align}
 \{\gamma_i,\gamma_j\}=\delta_{ij}
\end{align}
which are (up to a normalization) the relations of the Clifford algebra $\cliff_4$. These Clifford elements will be used as building blocks for a spinorial realization of $\os(4)$.
%

The Lie algebra $\os(4)$ is the algebra with $6$ generators, $\ell_{\mu\nu}$, $1\leq\mu<\nu\leq 4$, whose relations are given by
\begin{align}\label{eq_Lijrel}
 [\ell_{\mu\nu},\ell_{\rho\sigma}]=-i\left(\delta_{\nu\rho}\ell_{\mu\sigma}-\delta_{\nu\sigma}\ell_{\mu\rho}-\delta_{\mu\rho}\ell_{\nu\sigma}+\delta_{\mu\sigma}\ell_{\nu\rho}\right).
\end{align}
Let us denote $L_{\mu\nu}=a_\mu^{\dagger}a_\nu-a_\mu a_\nu^{\dagger}$ and $\Sigma_{\mu\nu}=\frac{1}{2}\gamma_\mu\gamma_\nu$. Both the combinations $-iL_{\mu\nu}$ and $-i\Sigma_{\mu\nu}$ realize the $\os(4)$ algebra. We now define the total angular momentum as the sum:
\begin{align}
 J_{\mu\nu}=-i(L_{\mu\nu}+\Sigma_{\mu\nu}).
\end{align}
These total angular momenta $J_{\mu\nu}$ realize again the $\os(4)$ commutation relations.

\subsection{The commutant}
We look for the commutant of the $\os(2)\oplus\os(2)$ subalgebra of $\os(4)$, that is, operators that commute with
\begin{align}\label{eq:J12J34}
\begin{aligned}{}
 J_{12}&=-i\left(a_{1}^{\dagger}a_2^{\phantom{\dagger}}-a_1^{\phantom{\dagger}}a_{2}^{\dagger}+\tfrac{1}{2}(b_1^{\phantom{\dagger}}+b_1^{\dagger})(b_2^{\phantom{\dagger}}+b_2^{\dagger})\right),\\
 J_{34}&=-i\left(a_{3}^{\dagger}a_4^{\phantom{\dagger}}-a_3^{\phantom{\dagger}}a_{4}^{\dagger}+\tfrac{1}{2}(b_3^{\phantom{\dagger}}+b_3^{\dagger})(b_4^{\phantom{\dagger}}+b_4^{\dagger})\right).
\end{aligned}
\end{align}
The combinations
\begin{align}
\begin{aligned}
 \mK_1&=\tfrac{1}{2}(a_1^{\dagger}a_1^{\phantom{\dagger}}+a_2^{\dagger}a_2^{\phantom{\dagger}}-a_3^{\dagger}a_3^{\phantom{\dagger}}-a_4^{\dagger}a_4^{\phantom{\dagger}}),\\
 \mK_2&=2(L_{12}\Sigma_{12}+L_{13}\Sigma_{13}+L_{14}\Sigma_{14}+L_{23}\Sigma_{23}+L_{24}\Sigma_{24}+L_{34}\Sigma_{34}-\tfrac{3}{4}),\\
 r&=i(b_1^{\phantom{\dagger}}+b_1^{\dagger})(b_2^{\phantom{\dagger}}+b_2^{\dagger}),\\
 \mR&=-(b_1^{\phantom{\dagger}}+b_1^{\dagger})(b_2^{\phantom{\dagger}}+b_2^{\dagger})(b_3^{\phantom{\dagger}}+b_3^{\dagger})(b_4^{\phantom{\dagger}}+b_4^{\dagger})
\end{aligned}
\end{align}
commute with $H$ and $J_{12}$, $J_{34}$. This is verified by a direct calculation.
The algebra generated by the elements $K_1$, $K_2$, $r$ closes onto the following form:
\begin{align}\label{eq:algebracommutant}
\begin{aligned}
 &[\mK_1,\mK_2]=\mK_3,\qquad\qquad [\mK_1,\mK_3]=\mK_2-(J_{12}+J_{34}\mR)r+\tfrac{1}{2},\\
 &[\mK_2,\mK_3]=4\mK_1[(J_{12}+J_{34}\mR)r-1]-2\mK_3(J_{12}+J_{34}\mR)r-2H(J_{12}-J_{34}\mR)r,\\
 &[\mK_1,r]=0,\qquad\qquad \{\mK_2,r\}=-r+2(J_{12}+J_{34}\mR)r,\qquad\qquad \{\mK_3,r\}=0.
\end{aligned}
\end{align}
Here $H$, $J_{12}$, $J_{34}$ and $\mR$ are central. These relations are identified with those in \eqref{eq:dualm1Hahn}. Hence we have obtained the centrally extended dual $-1$ Hahn algebra, or equivalently, the Schwinger-Dunkl algebra $\mathfrak{sd}(2)$, as a commutant.

\section{The algebra of the $\mathfrak{osp}(1|2)$ Clebsch-Gordan coefficients}\label{sec:clebschgordan}
In this Section, we first introduce the $\osp(1|2)$ algebra, present its Clebsch-Gordan problem, and then show how it is connected to the dual $-1$ Hahn algebra.

\subsection{The Lie superalgebra $\osp(1|2)$}
The $\osp(1|2)$ algebra can be presented as the algebra with generators $A_0$, $A_\pm$ and an involution $P$ encoding the $\mathbb{Z}_2$-grading of the superalgebra ($P$ commutes with the even element $A_0$ and anticommutes with the odd elements $A_\pm$). The defining relations are
\begin{align}\label{eq:osp12p1}
 \{A_+,A_-\}=2A_0,\qquad [A_0,A_\pm]=\pm A_\pm,\qquad  [P,A_0]=0,\qquad \{P,A_\pm\}=0.
\end{align}
The algebra $\osp(1|2)$ also possesses an sCasimir \cite{Lesniewski1995}
\begin{align}
 S=\tfrac{1}{2}\left([A_+,A_-]+1\right)=A_+A_--A_0+\tfrac{1}{2}
\end{align}
which commutes with the even elements and anticommutes with the odd elements
\begin{align}\label{eq:osp12p2}
 [S,A_0]=\{S,A_{\pm}\}=0.
\end{align}
Multiplying the sCasimir by the involution, we obtain a Casimir element for $\osp(1|2)$
\begin{align}
 Q=\left(A_+A_--A_0+\tfrac{1}{2}\right)P.
\end{align}
This Casimir element commutes with all generators of $\osp(1|2)$.

Positive infinite-dimensional discrete series representations for $\osp(1|2)$ are labelled by $(\mu,\epsilon)$, with $\mu\geq0$, $\epsilon=\pm1$. Let us denote by $\ket{n,\mu,\epsilon}$ the basis vectors associated to an irrep $(\mu,\epsilon)$. The generators act as follows in this basis:
\begin{subequations}\label{eq:actionosp12}
\begin{align}
\begin{aligned}
 A_0\ket{n,\mu,\epsilon}&=(n+\mu+\tfrac{1}{2})\ket{n,\mu,\epsilon},\\[.6em]
 P\ket{n,\mu,\epsilon}&=\epsilon(-1)^{n}\ket{n,\mu,\epsilon},
\end{aligned}\qquad\qquad
\begin{aligned}
 A_+\ket{n,\mu,\epsilon}&=\sqrt{[n+1]_\mu}\ket{n+1,\mu,\epsilon},\\
 A_-\ket{n,\mu,\epsilon}&=\sqrt{[n]_\mu}\ket{n-1,\mu,\epsilon},
\end{aligned}
\end{align}
where we define the \textit{mu}-numbers $[n]_\mu$ as
\begin{align}
 [n]_\mu=n+\mu\left(1-(-1)^{n}\right).
\end{align}
By Schur's lemma, the Casimir element acts as a multiple of the identity on these irreps
\begin{align}
 Q\ket{n,\mu,\epsilon}=-\epsilon\mu\ket{n,\mu,\epsilon}.
\end{align}
\end{subequations}

\subsection{The Clebsch-Gordan problem of $\osp(1|2)$}
We shall now look at the recoupling of two representations $(\mu_1,\epsilon_1)$ and $(\mu_2,\epsilon_2)$ of $\osp(1|2)$.

The direct product representation $(\mu_1,\epsilon_1)\otimes(\mu_2,\epsilon_2)$ has the associated basis vectors\!\!\!  \penalty-10000$\ket{n_1,\mu_1,\epsilon_1}\otimes\ket{n_2,\mu_2,\epsilon_2}$. 
The elements $A_0\otimes1$, $1\otimes A_0$, $P\otimes1$, $1\otimes P$ as well as the Casimir elements $Q\otimes1$, $1\otimes Q$ are diagonal in this basis. Equivalently, we shall consider the elements $(A_0\otimes1+1\otimes A_0)$, $(A_0\otimes1-1\otimes A_0)$, $P\otimes1$ and $P\otimes P$, also diagonal in this basis.

This direct product representation admits the following decomposition \cite{Genest2015a} in irreducibles:
\begin{align}\label{eq:irrepsdecomposition}
 (\mu_1,\epsilon_1)\otimes(\mu_2,\epsilon_2)=\bigoplus_{j=0}^{\infty}(\mu_{12},\epsilon_{12}),
\end{align}
where
\begin{align}
 \mu_{12}=\mu_1+\mu_2+j+\tfrac{1}{2},\qquad \epsilon_{12}=(-1)^{j}\epsilon_1\epsilon_2.
\end{align}
The coupled basis vectors associated to the irreducibles $(\mu_{12},\epsilon_{12})$ are denoted $\ket{n_{12},\mu_{12},\epsilon_{12}}$.
The Casimir elements $Q\otimes1$, $1\otimes Q$ are again diagonal in the coupled basis. The other diagonal elements can be obtained as follows:

Consider the coproduct map $\Delta:\osp(1|2)\to\osp(1|2)\otimes\osp(1|2)$, which is a coassociative algebra morphism. The coproduct maps the $\osp(1|2)$ generators as follows
\begin{align}
\begin{aligned}
 \Delta(A_0)=A_0^{(12)}&=A_0\otimes1+1\otimes A_0\\
 \Delta(A_\pm)=A_\pm^{(12)}&=A_\pm\otimes P+1\otimes A_\pm\\
 \Delta(P)=P^{(12)}&=P\otimes P
\end{aligned}\!\!
\begin{aligned}
 &=A_0^{(1)}+A_0^{(2)},\\
 &=A_\pm^{(1)}P^{(2)}+A_\pm^{(2)},\\
 &=P^{(1)}P^{(2)},
\end{aligned}
\end{align}
and the Casimir element according to
\begin{align}\label{eq:Q12algebra}
 \Delta(Q)=Q^{(12)}=(A_-^{(1)}A_+^{(2)}-A_+^{(1)}A_-^{(2)})P^{(1)}+Q^{(1)}P^{(2)}+Q^{(2)}P^{(1)}-\tfrac{1}{2}P^{(1)}P^{(2)},
\end{align}
where the superindex denotes on which factor of $\osp(1|2)$ the operator is acting. The elements $\Delta(A_0)=A_0\otimes1+1\otimes A_0$, $\Delta(P)=P\otimes P$ and $\Delta(Q)$ are diagonal in the coupled basis.
\\[1em]
The decomposition \eqref{eq:irrepsdecomposition} indicates that the vector spaces spanned by $\ket{n_1,\mu_1,\epsilon_1}\otimes\ket{n_2,\mu_2,\epsilon_2}$ and $\ket{n_{12},\mu_{12},\epsilon_{12}}$ are isomorphic; one defines the Clebsch-Gordan coefficients $\mathcal{C}_{n_{12},j}^{n_1,n_2}$ as the expansion coefficients between the two bases
\begin{align}
 \ket{n_{12},\mu_{12},\epsilon_{12}}=\sum_{n_1,n_2}\mathcal{C}_{n_{12},j}^{n_1,n_2}\ket{n_1,\mu_1,\epsilon_1}\otimes\ket{n_2,\mu_2,\epsilon_2}.
\end{align}
The Clebsch-Gordan coefficients are characterized algebraically by the elements that are diagonalized in each of the bases. In particular, $P\otimes1$, $(A_0\otimes1-1\otimes A_0)$ and $\Delta(Q)$ are not diagonal in both bases, so they obey non-trivial commutation relations. The algebra formed by these elements determines the Clebsch-Gordan coefficients.

Write
\begin{align}
 \kappa_1=\tfrac{1}{2}(A_0\otimes1-1\otimes A_0),\qquad \kappa_2=Q^{(12)}P^{(12)},\qquad \text{and}\quad p=P\otimes1,
\end{align}
a straightforward calculation in $\osp(1|2)\otimes\osp(1|2)$ yields
\begin{align}
\begin{aligned}
 &[\kappa_1,\kappa_2]=\kappa_3,\qquad\qquad [\kappa_1,\kappa_3]=\kappa_2-(Q^{(1)}+Q^{(2)}P^{(12)})p+\tfrac{1}{2},\\
 &[\kappa_3,\kappa_2]=4\kappa_1\!\!\left(1-(Q^{(1)}\!+\!Q^{(2)}P^{(12)})p\right)\!+2p\!\left(\kappa_3(Q^{(1)}\!+\!Q^{(2)}P^{(12)})+A_0^{(12)}(Q^{(1)}\!-\!Q^{(2)}P^{(12)})\right)\!,\\
 &[\kappa_1,p]=0,\qquad\qquad \{\kappa_2,p\}=-p+2(Q^{(1)}\!+\!Q^{(2)}P^{(12)}),\qquad\qquad \{\kappa_3,p\}=0.
\end{aligned}
\end{align}
Keeping in mind that the elements $Q^{(1)}$, $Q^{(2)}$, $P^{(12)}$ and $A_0^{(12)}$ are central since they are diagonalized in both bases, one recognizes the defining relations of the (centrally extended) dual $-1$ Hahn algebra \eqref{eq:dualm1Hahn}. This reveals that the Clebsch-Gordan coefficients of $\osp(1|2)$ are (essentially) the dual $-1$ Hahn polynomials \cite{Tsujimoto2011,Genest2013,Bergeron2016}.

\section{The Howe duality correspondence}\label{sec:howeduality}

In Sections \ref{sec:dualcommutant} and \ref{sec:clebschgordan}, we have obtained the (centrally extended) dual $-1$ Hahn algebra in two different contexts. We will now reinterpret the contents of the last two Sections in order to display the two presentations of the algebra in a unified way.

\subsection{Connecting the two approaches}
Let us go back to our construction in Section \ref{sec:dualcommutant} involving four Bosonic and Fermionic oscillators. We can form three copies of $\osp(1|2)$ labelled by one of the sets $S\in\{\{1,2\},\{3,4\},\{1,2,3,4\}\}$ by introducing the operators:
\begin{align}\label{eq:osp12howe}
 \mA_-^{S}=\sum_{\mu\in S}a_\mu\gamma_\mu,\qquad \mA_+^{S}=\sum_{\mu\in S}a_\mu^{\dagger}\gamma_\mu,\qquad \mA_0^{S}=\frac{1}{2}\sum_{\mu\in S}\{a_\mu^{\dagger},a_\mu^{\phantom{\dagger}}\},\qquad \mP^{S}=e^{i\pi|S|/4}\prod_{\mu\in S}\gamma_\mu,
\end{align}
which obey the defining relations of $\osp(1|2)$ given in \eqref{eq:osp12p1}, \eqref{eq:osp12p2}. The Casimir element associated to each set $S$ is given by
\begin{align}
 \mQ^{S}=(\mA_+^{S}\mA_-^{S}-\mA_0^{S}+\tfrac{1}{2})\mP^{S}.
\end{align}
Let us revisit the Clebsch-Gordan problem of $\osp(1|2)$ in this framework. As seen in Section \ref{sec:clebschgordan}, the elements of the set
\begin{align}
 \mathcal{E}=\left\{(\mA_0^{\{1,2\}}+\mA_0^{\{3,4\}}), (\mA_0^{\{1,2\}}-\mA_0^{\{3,4\}}), \mP^{\{1,2\}}, \mP^{\{1,2\}}\mP^{\{3,4\}}, \mQ^{\{1,2\}}, \mQ^{\{3,4\}}, \mQ^{\{1,2,3,4\}}\right\}
\end{align}
are diagonal in at least one of the bases. The elements that are diagonal in both bases commute with all others, and the remaining ones obey the commutation relations of the dual $-1$ Hahn algebra.

We now give the explicit form of all elements in $\mathcal{E}$ in the realization \eqref{eq:osp12howe}. It will appear that they can be matched with certain expressions given in Section \ref{sec:dualcommutant}, thus explaining why the dual $-1$ Hahn algebra appeared in two seemingly different situations. The expressions are the following once translated in terms of the $a$ and $b$ ladder operators:
\begin{align}\label{eq:Dduals}
\begin{aligned}
 \mA_0^{\{1,2\}}+\mA_0^{\{3,4\}}&=H,\\
 \mA_0^{\{1,2\}}-\mA_0^{\{3,4\}}&=2K_1,\\
 \mP^{\{1,2\}}&=r,\\
 \mP^{\{1,2\}}\mP^{\{3,4\}}&=P_{1234},\\
 \mQ^{\{1,2\}}&=J_{12},\\
 \mQ^{\{3,4\}}&=J_{34},\\
 \mQ^{\{1,2,3,4\}}&=K_2P_{1234}.
\end{aligned}
\end{align}
In this framework, we easily see that the statement that $\{\mQ^{\{1,2\}},\mQ^{\{3,4\}}\}$ commute with all other generators implies that the commutant of $\{J_{12}, J_{34}\}$ contains the dual $-1$ Hahn algebra.

\subsection{An instance of Howe duality}
That the dual $-1$ Hahn algebra can be viewed on the one hand as the commutant of $\mathfrak{o}(2) \oplus \mathfrak{o}(2)$ in a spinorial representation of $\mathfrak{o}(4)$ and can be embedded in $U(\mathfrak{osp}(1|2)) \otimes U(\mathfrak{osp}(1|2))$ on the other hand can be attributed to Howe duality as we now explain.

It is known \cite{Brackx2008} that $\osp(1|2)$ and $Pin(2n)$ have dual (commuting) actions on the space of polynomials $\mathcal{P}(\mathbb{R}^{m},\mathbb{S})$ defined in Euclidean space $\mathbb{R}^{m}$ and taking values in a spinor space $\mathbb{S}$. In our situation, at the level of the algebras, this would correspond to the fact that the generators of $\osp(1|2)$ and $\os(2n)$ commute. A direct computation using the expressions \eqref{eq:J12J34} and \eqref{eq:osp12howe} confirms that it is indeed the case:
\begin{align}
 &[J_{12},\mA_\bullet^{\{1,2\}}]=[J_{34},\mA_\bullet^{\{3,4\}}]=0,\\
 [J_{ij},\,&\mA_\bullet^{\{1,2,3,4\}}]=0,\qquad 1\leq i<j\leq4.
\end{align}
As a byproduct, it can be shown that the Casimir elements of both algebras can be put in correspondance. Recall that the Casimir element of $\os(2n)$ denoted $C^{\{1,...,2n\}}$ is given by
\begin{align}
 C^{\{1,...,2n\}}=\sum_{1\leq i<j\leq2n}{J_{ij}}^{2}.
\end{align}
The Casimir elements of both $\osp(1|2)$ and $\os(2n)$ associated to each of the three copies labelled by $S$ are related by
\begin{align}
 C^{\{1,2\}}=(\mQ^{\{1,2\}})^{2},\qquad C^{\{3,4\}}=(\mQ^{\{3,4\}})^{2},\qquad C^{\{1,2,3,4\}}=(\mQ^{\{1,2,3,4\}})^{2}-\tfrac{3}{4}.
\end{align}
Remark: That the relation between the Casimir element of $\os(2n)$ and the one of $\osp(1|2)$ is quadratic does not come as a surprise. It is known \cite{Frappat2019a} that in the context of Howe duality for the pair $(\su(1,1),\os(2n))$, the relation between the associated Casimir elements is linear, whilst the algebra $\osp(1|2)$ can be seen as a ``square root'' of $\su(1,1)$.

\section{Conclusion}\label{sec:conclusion}
To sum up, we have presented two frameworks that lead to the (centrally extended) dual $-1$ Hahn algebra: one in which we looked at the commutant of a spinorial realization of the $\os(2)\oplus\os(2)$ subalgebra of $\os(4)$, and the other in which we looked at the algebra of the $\osp(1|2)$ Clebsch-Gordan coefficients. We have explained how these two approaches are dual (in the sense of Howe) by considering representations that featured the dual pair $(\osp(1|2),\os(2n))$. We have also highlighted how the results presented in this report can be seen as a ``square root'' of those related to the dual pair $(\su(1,1),\os(2n))$ \cite{Frappat2019a}.

One should note that the construction presented here generalizes straightforwardly if one considers instead the commutant of the spinorial representation of the $\os(m)\oplus\os(m')$ subalgebra of $\os(m+m')$; the (centrally extended) dual $-1$ Hahn algebra is still recovered and the Howe duality again operates.

\subsection*{Acknowledgments}
The authors benefitted from discussions with Nicolas Crampé, Hendrik De Bie, Luc Frappat, Eric Ragoucy, Stéphane Vinet and Alexei Zhedanov.
JG holds an Alexander-Graham-Bell scholarship from the Natural Science and Engineering Research Council (NSERC) of Canada.
The research of LV is supported in part by a Discovery Grant from NSERC.

\bibliographystyle{unsrtinurl} 
\bibliography{dual-1Hcomm.bib} 

\begin{thebibliography}{10}

\bibitem{Genest2013}
V.~X. Genest, L.~Vinet, and A.~Zhedanov.
\newblock {The algebra of dual {$-1$} Hahn polynomials and the Clebsch–Gordan
  problem of {$\mathfrak{sl}_{-1}(2)$}}.
\newblock {\em Journal of Mathematical Physics}, 54(2):1--15, 2013.
\newblock \href {http://arxiv.org/abs/1207.4220} {\path{arXiv:1207.4220}}.

\bibitem{Tsujimoto2013}
S.~Tsujimoto, L.~Vinet, and A.~Zhedanov.
\newblock {Dual {$-1$} Hahn polynomials: “Classical” polynomials beyond the
  Leonard duality}.
\newblock {\em Proceedings of the American Mathematical Society},
  141(3):959--970, 2013.

\bibitem{Koekoek2010}
R.~Koekoek, P.~A. Lesky, and R.~F. Swarttouw.
\newblock {\em {Hypergeometric Orthogonal Polynomials and Their
  {$q$}-Analogues}}.
\newblock Springer Monographs in Mathematics. Springer Berlin Heidelberg, 2010.

\bibitem{Tsujimoto2012}
S.~Tsujimoto, L.~Vinet, and A.~Zhedanov.
\newblock {Dunkl shift operators and Bannai-Ito polynomials}.
\newblock {\em Advances in Mathematics}, 229(4):2123--2158, 2012.
\newblock \href {http://arxiv.org/abs/1106.3512} {\path{arXiv:1106.3512}}.

\bibitem{DeBie2016a}
H.~{De Bie}, V.~X. Genest, and L.~Vinet.
\newblock {The {$\mathbb{Z}_2^n$} Dirac–Dunkl operator and a higher rank
  Bannai–Ito algebra}.
\newblock {\em Advances in Mathematics}, 303:390--414, 2016.
\newblock \href {http://arxiv.org/abs/1511.02177} {\path{arXiv:1511.02177}}.

\bibitem{DeBie2017a}
H.~{De Bie}, V.~X. Genest, W.~van~de Vijver, and L.~Vinet.
\newblock {A higher rank Racah algebra and the {$\mathbb{Z}_2^n$}
  Laplace–Dunkl operator}.
\newblock {\em Journal of Physics A: Mathematical and Theoretical},
  51(2):025203, 2017.
\newblock \href {http://arxiv.org/abs/1610.02638} {\path{arXiv:1610.02638}}.

\bibitem{Post2017}
S.~Post and A.~Walter.
\newblock {A higher rank extension of the Askey-Wilson Algebra}.
\newblock 2017.
\newblock \href {http://arxiv.org/abs/1705.01860v2}
  {\path{arXiv:1705.01860v2}}.

\bibitem{DeBie2018}
H.~{De Bie}, H.~{De Clercq}, and W.~{Van De Vijver}.
\newblock {The higher rank {$q$}-deformed Bannai-Ito and Askey-Wilson algebra}.
\newblock 2018.
\newblock \href {http://arxiv.org/abs/1805.06642v1}
  {\path{arXiv:1805.06642v1}}.

\bibitem{Crampe2019}
N.~Cramp{\'{e}}, E.~Ragoucy, L.~Vinet, and A.~Zhedanov.
\newblock {Truncation of the reflection algebra and the Hahn algebra}.
\newblock {\em Journal of Physics A: Mathematical and Theoretical},
  52(35):1--9, 2019.
\newblock \href {http://arxiv.org/abs/1903.05674} {\path{arXiv:1903.05674}}.

\bibitem{Crampe2019d}
N.~Cramp{\'{e}}, L.~P. D'Andecy, and L.~Vinet.
\newblock {Temperley-Lieb, Brauer and Racah algebras and other centralizers of
  {$\mathfrak{su}(2)$}}.
\newblock 2019.
\newblock \href {http://arxiv.org/abs/1905.06346} {\path{arXiv:1905.06346}}.

\bibitem{Crampe2020}
N.~Cramp{\'{e}}, J.~Gaboriaud, L.~Vinet, and M.~Zaimi.
\newblock {Revisiting the Askey–Wilson algebra with the universal
  {$R$}-matrix of {$U_q(\mathfrak{sl}_2)$}}.
\newblock {\em Journal of Physics A : Mathematical and Theoretical}, 23:05LT01,
  2020.
\newblock \href {http://arxiv.org/abs/1908.04806} {\path{arXiv:1908.04806}}.

\bibitem{Terwilliger2014}
P.~Terwilliger and A.~{\v{Z}}itnik.
\newblock {Distance-regular graphs of {$q$}-Racah type and the universal
  Askey-Wilson algebra}.
\newblock {\em Journal of Combinatorial Theory. Series A}, 125(1):98--112,
  2014.
\newblock \href {http://arxiv.org/abs/1307.7968} {\path{arXiv:1307.7968}}.

\bibitem{Bullock1999}
D.~Bullock and J.~H. Przytycki.
\newblock {Multiplicative structure of Kauffman bracket skein module
  quantizations}.
\newblock {\em Proceedings of the American Mathematical Society},
  128(3):923--931, 1999.
\newblock \href {http://arxiv.org/abs/math/9902117}
  {\path{arXiv:math/9902117}}.

\bibitem{Genest2014}
V.~X. Genest, L.~Vinet, and A.~Zhedanov.
\newblock {The Racah algebra and superintegrable models}.
\newblock {\em Journal of Physics: Conference Series}, 512(1):012011, 2014.
\newblock \href {http://arxiv.org/abs/1312.3874} {\path{arXiv:1312.3874}}.

\bibitem{DeBie2015}
H.~{De Bie}, V.~X. Genest, S.~Tsujimoto, L.~Vinet, and A.~Zhedanov.
\newblock {The Bannai-Ito algebra and some applications}.
\newblock {\em Journal of Physics: Conference Series}, 597:012001, 2015.
\newblock \href {http://arxiv.org/abs/1411.3913} {\path{arXiv:1411.3913}}.

\bibitem{Baseilhac2005}
P.~Baseilhac.
\newblock {Deformed Dolan-Grady relations in quantum integrable models}.
\newblock {\em Nuclear Physics B}, 709(3):491--521, 2005.
\newblock \href {http://arxiv.org/abs/hep-th/0404149v3}
  {\path{arXiv:hep-th/0404149v3}}.

\bibitem{Gaboriaud2018}
J.~Gaboriaud, L.~Vinet, S.~Vinet, and A.~Zhedanov.
\newblock {The Racah algebra as a commutant and Howe duality}.
\newblock {\em Journal of Physics A : Mathematical and Theoretical}, 51:50LT01,
  2018.
\newblock \href {http://arxiv.org/abs/1808.05261} {\path{arXiv:1808.05261}}.

\bibitem{Gaboriaud2018a}
J.~Gaboriaud, L.~Vinet, S.~Vinet, and A.~Zhedanov.
\newblock {The generalized Racah algebra as a commutant}.
\newblock {\em Journal of Physics: Conference Series}, 1194:012034, 2019.
\newblock \href {http://arxiv.org/abs/1808.09518v1}
  {\path{arXiv:1808.09518v1}}.

\bibitem{Gaboriaud2018b}
J.~Gaboriaud, L.~Vinet, S.~Vinet, and A.~Zhedanov.
\newblock {The dual pair {$Pin(2n)\otimes\mathfrak{osp}(1|2)$}, the Dirac
  equation and the Bannai–Ito algebra}.
\newblock {\em Nuclear Physics B}, 937:226--239, 2018.
\newblock \href {http://arxiv.org/abs/1810.00130} {\path{arXiv:1810.00130}}.

\bibitem{Frappat2019}
L.~Frappat, J.~Gaboriaud, E.~Ragoucy, and L.~Vinet.
\newblock {The {$q$}-Higgs and Askey-Wilson algebras}.
\newblock {\em Nuclear Physics B}, 1:114632, 2019.
\newblock \href {http://arxiv.org/abs/1903.04616} {\path{arXiv:1903.04616}}.

\bibitem{Gaboriaud2019}
J.~Gaboriaud, L.~Vinet, and S.~Vinet.
\newblock {Howe duality and algebras of the Askey-Wilson type: an overview}.
\newblock (6):1--8, 2019.
\newblock \href {http://arxiv.org/abs/1911.08314} {\path{arXiv:1911.08314}}.

\bibitem{Higgs1979}
P.~W. Higgs.
\newblock {Dynamical symmetries in a spherical geometry I}.
\newblock {\em Journal of Physics A: Mathematical and General}, 12(3):309,
  1979.

\bibitem{Frappat2019a}
L.~Frappat, J.~Gaboriaud, L.~Vinet, S.~Vinet, and A.~Zhedanov.
\newblock {The Higgs and Hahn algebras from a Howe duality perspective}.
\newblock {\em Physics Letters, Section A: General, Atomic and Solid State
  Physics}, 383(14):1531--1535, 2019.
\newblock \href {http://arxiv.org/abs/1811.09359} {\path{arXiv:1811.09359}}.

\bibitem{Genest2012}
V.~X. Genest, M.~E. Ismail, L.~Vinet, and A.~Zhedanov.
\newblock {The Dunkl oscillator in the plane: I. Superintegrability, separated
  wavefunctions and overlap coefficients}.
\newblock {\em Journal of Physics A: Mathematical and Theoretical}, 46(14),
  2013.
\newblock \href {http://arxiv.org/abs/1212.4459} {\path{arXiv:1212.4459}}.

\bibitem{Genest2014a}
V.~X. Genest, M.~E.~H. Ismail, L.~Vinet, and A.~Zhedanov.
\newblock {The Dunkl Oscillator in the Plane II: Representations of the
  Symmetry Algebra}.
\newblock {\em Communications in Mathematical Physics}, 329(3):999--1029, 2014.
\newblock \href {http://arxiv.org/abs/1302.6142} {\path{arXiv:1302.6142}}.

\bibitem{Bernard2020}
P.-A. Bernard, J.~Gaboriaud, and L.~Vinet.
\newblock {Superintegrability and the dual {$-1$} Hahn algebra in
  superconformal quantum mechanics}.
\newblock pages 1--24, 2020.
\newblock \href {http://arxiv.org/abs/2001.07309} {\path{arXiv:2001.07309}}.

\bibitem{Lesniewski1995}
A.~Le{\'{s}}niewski.
\newblock {A remark on the Casimir elements of Lie superalgebras and quantized
  Lie superalgebras}.
\newblock {\em Journal of Mathematical Physics}, 36(3):1457--1461, 1995.

\bibitem{Genest2015a}
V.~X. Genest, L.~Vinet, and A.~Zhedanov.
\newblock {A Laplace-Dunkl Equation on {$S^2$} and the Bannai–Ito Algebra}.
\newblock {\em Communications in Mathematical Physics}, 336(1):243--259, 2015.
\newblock \href {http://arxiv.org/abs/1312.6604} {\path{arXiv:1312.6604}}.

\bibitem{Tsujimoto2011}
S.~Tsujimoto, L.~Vinet, and A.~Zhedanov.
\newblock {From {$sl_q(2)$} to a Parabosonic Hopf Algebra}.
\newblock {\em Symmetry, Integrability and Geometry: Methods and Applications},
  7(093):13, 2011.
\newblock \href {http://arxiv.org/abs/1108.1603v3} {\path{arXiv:1108.1603v3}}.

\bibitem{Bergeron2016}
G.~Bergeron and L.~Vinet.
\newblock {Generating functions for the {$\mathfrak{osp}(1|2)$}
  Clebsch–Gordan coefficients}.
\newblock {\em Journal of Physics A: Mathematical and Theoretical},
  49(11):115202, 2016.
\newblock \href {http://arxiv.org/abs/1507.00018} {\path{arXiv:1507.00018}}.

\bibitem{Brackx2008}
F.~Brackx, H.~{De Schepper}, D.~Eelbode, and V.~Sou{\v{c}}ek.
\newblock {The Howe Dual Pair in Hermitean Clifford Analysis}.
\newblock {\em Revista Matem{\'{a}}tica Iberoamericana}, 26(2):449--479, 2010.

\end{thebibliography}

\end{document}